\documentclass[11pt]{article}

\usepackage{amssymb,graphicx}

\def\emline#1#2#3#4#5#6{%
       \put(#1,#2){\special{em:moveto}}%
       \put(#4,#5){\special{em:lineto}}}
\def\newpic#1{}

\newtheorem{dfn}{Definition}
\newtheorem{thm}[dfn]{Theorem}
\newtheorem{lem}[dfn]{Lemma}

\newtheorem{conj}[dfn]{Conjecture}
\newtheorem{athm}{Theorem}

\newtheorem{thmm}{Theorem}

\newtheorem{prerem}{Remark}

\newtheorem{preproof}{{\bf Proof}}

\newenvironment{proof}[1]{\begin{preproof}{\rm #1}
                \hfill{$\blacksquare$}}{\end{preproof}}
\newtheorem{preprooff}{{\bf Proof}}

\newenvironment{prooff}[1]{\begin{preprooff}{\rm #1}}{\end{preprooff}}
\newcommand{\utlc}{uniquely $2$-list colorable}
\newcommand{\uthreelc}{uniquely $3$-list colorable}
\newcommand{\uklc}{uniquely $k$-list colorable}

\newcommand{\uttlc}{uniquely $(2,t)$-list colorable}
\newcommand{\uktlc}{uniquely $(k,t)$-list colorable}
\newcommand{\uftlc}{uniquely $(f,t)$-list colorable}
\newcommand{\deffnt}{\sf}
\newcommand{\reffnt}{\rm}
\newcommand{\N}{\Bbb{N}}
\def\Chi{\lower-.3ex\hbox{$\chi$}}
\def\le{\leqslant}
\def\ge{\geqslant}
\title{Uniquely $2$--List Colorable Graphs\footnote{
The research of second and third authors is supported by the
Institute for Studies in Theoretical Physics and Mathematics (IPM),
Tehran, Iran.}}
\author{ 
\bf Y.~Ganjali~G., M.~Ghebleh, H.~Hajiabolhassan,
\vspace{2mm}\\ %
\bf M.~Mirzazadeh, and B.\,S.~Sadjad
\vspace{3eX}\\ %
\small Institute for Studies in Theoretical Physics \\ %
\small and Mathematics (IPM), Tehran, Iran %
\vspace{1eX}\\ %
\small and
\vspace{1eX}\\ %
\small Department of Mathematical Sciences\\ %
\small Sharif University of Technology\\ %
\small P.\,O.~Box~11365--9415, Tehran, Iran\\ %
}
\date{}
\begin{document}
\maketitle

\begin{abstract}
\addtolength{\baselineskip}{.7mm}
A graph is called to be uniquely list colorable, if it admits a
list assignment which induces a unique list coloring. We study
uniquely list colorable graphs with a restriction on the number
of colors used. In this way we generalize a theorem which
characterizes uniquely $2$--list colorable graphs. We introduce
the uniquely list chromatic number of a graph and make a
conjecture about it which is a generalization of the well known
Brooks' theorem.

\addtolength{\baselineskip}{-.7mm}
\end{abstract}
\addtolength{\baselineskip}{.8mm}
\section{Introduction}
We consider finite, undirected simple graphs. For necessary
definitions and notations we refer the reader to standard texts
such as \cite{west}.

Let $G$ be a graph, $f:V(G)\to\N$ be a given map, and $t\in\N$.
An {\deffnt $(f,t)$-list assignment} $L$ to $G$
is a map, which assigns to each vertex $v$, a set $L(v)$ of size
$f(v)$ and $|\bigcup_vL(v)|=t$. By a {\deffnt list coloring} for
$G$ from such $L$ or an {\deffnt $L$-coloring} for short, we
shall mean a proper coloring $c$ in which $c(v)$ is chosen from
$L(v)$, for each vertex $v$. When $f(v)=k$ for all $v$, we simply
say $(k,t)$-list assignment for an $(f,t)$-list assignment.
When the parameter $t$ is not of special interest, we say
$f$-list (or $k$-list) assignment simply. Specially if $L$ is
a $(t,t)$-list assignment to $G$, then any $L$-coloring is
called a $t$-coloring for~$G$.

In this paper we study the concept of
uniquely list coloring which was introduced by
Dinitz and Martin~\cite{dinmar} and independently by
Mahdian and Mahmoodian~\cite{mahmah}. In~\cite{dinmar} and
\cite{mahmah} uniquely $k$-list colorable 
graphs are
introduced as graphs who admit a
$k$-list assignment which induces a unique list coloring. In the
present work we study uniquely list colorings of graphs in a more
general sense.

\begin{dfn}
Suppose that $G$ is a graph, $f:V(G)\to\N$ is a map, and $t\in\N$.
The graph $G$ is called to be {\deffnt uniquely $(f,t)$-list
colorable} if there exists an $(f,t)$-list assignment $L$ to $G$,
such that $G$ has a unique $L$-coloring. We call $G$ to be
{\deffnt uniquely $f$-list colorable} if it is \uftlc\ for some~$t$.
\end{dfn}

If $G$ is a uniquely $(f,t)$-list (resp. $f$-list) colorable graph
and $f(v)=k$ for each $v\in V(G)$, we simply say that $G$ is a
uniquely $(k,t)$-list (resp. $k$-list) colorable graph.
In~\cite{mahmah} all \utlc\ graphs are characterized as follows.

\begin{athm}{\reffnt\cite{mahmah}}
\label{mahu2lc} A graph $G$ is not \utlc, if and only
if each of its blocks is either a complete graph, a complete
bipartite graph, or a cycle.
\end{athm}

For recent advances in uniquely list colorable graphs we direct
the interested reader to \cite{ghbmah} and \cite{egh}.

In developing computer programs for recognition of uniquely
$k$-list colorability of graphs, it is important to restrict the
number of colors as much as possible. So if $G$ is a \uklc\ graph,
the minimum number of colors which are sufficient for a $k$-list
assignment to $G$ with a unique list coloring,
will be an important parameter for us. Uniquely list colorable
graphs are related to defining sets of graph colorings as
discussed in~\cite{mahmah}, and in this application also the
number of colors is an important quantity.

In next section we show that for every \utlc\ graph $G$ there
exists a $2$-list assignment $L$, such that $G$
has a unique $L$-coloring and there are $\max\{3,\Chi(G)\}$
colors used in $L$.

\section{Uniquely $(2,t)$-list colorable graphs}

It is easy to see that for each \uklc\ graph $G$, and each
$k$-list assignment $L$ to its vertices which induces a unique
list coloring, at least $k+1$ colors must be used in $L$, and on
the other hand since $G$ has an $L$-coloring, at least $\Chi(G)$
colors must be used. So the number of colors used is at least
$\max\{k+1,\Chi(G)\}$ colors. Throughout this section our goal
is to prove the following theorem which implies the equality in
the case $k=2$.
\begin{thmm}
A graph $G$ is \utlc\ if and only if it is \uttlc, where
$t=\max\{3,\Chi(G)\}$.
\end{thmm}

To prove the theorem above we consider a counterexample $G$ to the
statement with minimum number of vertices. In theorems~\ref{block},
\ref{tfree}, and \ref{cyclechord}, we will show that $G$ is
$2$-connected and triangle-free, and each of its cycles is induced
(chordless).

As mentioned above, if $G$ is a \uklc\ graph,
and $L$ a $(k, t)$-list assignment to $G$ such
that $G$ has a unique $L$-coloring, then
$t\ge\max\{k+1,\Chi(G)\}$. Although the theorem above states that
when $k=2$ there exists an $L$ for which equality holds, this is not
the case in general.

\begin{figure}
\centering
\includegraphics{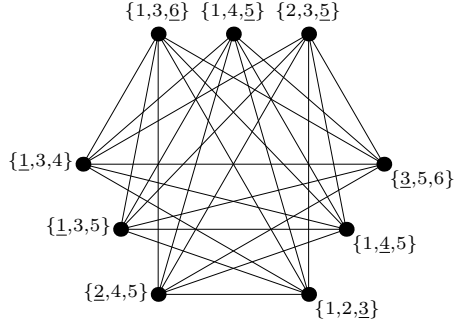}
\caption{A $3$-list assignment to $K_{3,3,3}$ which induces a unique
list coloring}
\label{K333}
\end{figure}

To see this, consider a complete tripartite \uthreelc\ graph $G$.
We will call each of the three color classes of $G$ a part.
In~\cite{ghbmah} it is shown that for each $k\ge 3$ there exists a
complete tripartite \uklc\ graph. For example one can check that
the graph $K_{3,3,3}$ has a unique list coloring from the lists
shown in Figure~\ref{K333} (the color taken by each vertex is
underlined).

Suppose that $L$ is a $(3,t)$-list assignment to $G$
which induces a unique list coloring~$c$, and the vertices of a
part $X$ of $G$ take on the same color $i$ in~$c$.
We introduce a $2$-list assignment $L'$ to $G\setminus X$ as follows.
For every vertex $v$ in $G\setminus X$, if $i\in L(v)$ then
$L'(v)=L(v)\setminus\{i\}$, and otherwise $L'(v)=L(v)\setminus\{j\}$
where $j\in L(v)$ and $j\not=c(v)$. Since $L$ induces a unique list
coloring $c$ for $G$, $G\setminus X$ has exactly one $L'$-coloring,
namely the restriction of $c$ to $V(G)\setminus X$. But $G\setminus X$
is a complete bipartite graph and this contradicts Theorem~\ref{mahu2lc}.
So on each part of $G$ there must be appeared at least $2$ colors
and therefore we have $t\ge 6$ while $\max\{k+1,\Chi(G)\}=4$.

Similarly one can see that if $G$ is a complete tripartite \uklc\
graph for some $k\ge 3$, and $L$ a $(k,t)$-list assignment to $G$
which induces a unique list coloring, then on each part there are at
least $k-1$ colors appeared and so we have $t\ge 3(k-1)$ while
$\max\{k+1,\Chi(G)\}=k+1$.

Towards our main theorem, we start with two basic lemmas.

\begin{lem}
\label{one1}Suppose that $G$ is a connected graph and
$f:V(G)\to\{1,2\}$ such that $f(v_0)=1$ for some vertex $v_0$ of $G$.
Then $G$ is a uniquely $(f,\Chi(G))$-list colorable graph.
\end{lem}

\begin{proof}{
Consider a spanning tree $T$ in $G$ rooted at $v_0$ and consider a
$\Chi(G)$-coloring $c$ for $G$. Let $L(v)$ be $\{c(v)\}$ if
$f(v)=1$, and $\{c(u),c(v)\}$ if $f(v)=2$, where $u$ is the parent
of $v$ in $T$. It is easy to see that $c$ is the only
$L$-coloring of $G$.
}\end{proof}

\begin{lem}
\label{cutvertex}Let $G$ be the union of two graphs $G_1$ and
$G_2$ which are joined in exactly one vertex $v_0$. Then $G$ is
\uttlc\ if and only if at least one of $G_1$ and $G_2$ is \uttlc.
\end{lem}

\begin{proof}{ %
If either $G_1$ or $G_2$ is a \uttlc\ graph, by use of
Lemma~\ref{one1} it is obvious that $G$ is also \uttlc. On the
other hand suppose that none of $G_1$ and $G_2$ is a \uttlc\
graph and $L$ is a $(2,t)$-list assignment to $G$
which induces a list coloring $c$. Since $G_1$ and $G_2$ are not
\uttlc, each of these has another coloring, say $c_1$ and $c_2$
respectively. If $c_1(v_0)=c(v_0)$ or $c_2(v_0)=c(v_0)$ then an
$L$-coloring for $G$ different from $c$ is obtained obviously.
Otherwise $c_1(v_0)=c_2(v_0)$, so we obtain a new $L$-coloring
for $G$, by combining $c_1$ and~$c_2$.
}\end{proof}

The following theorem is immediate by Lemma~\ref{one1} and
Lemma~\ref{cutvertex}.
\begin{thm}
Suppose that $G$ is a graph and $t\ge \Chi(G)$. The
graph $G$ is \uttlc\ if and only if at least one of its blocks is
a \uttlc\ graph.
\label{block}
\end{thm}

Next lemma which is an obvious statement, is useful throughout
the paper.
\begin{lem}
\label{adduv}Suppose that the independent vertices $u$ and $v$ in
a graph $G$ take on different colors in each $t$-coloring of $G$.
Then the graph $G$ is \uftlc\ if and only if $G+uv$ is a \uftlc\
graph.
\end{lem}

\begin{figure}
\centering %
\footnotesize %
\special{em:linewidth 0.4pt}
\unitlength 1.00mm
\linethickness{0.4pt}
\begin{picture}(47.00,37.00)
\put(4.33,21.00){\circle*{2.00}}
\put(44.33,21.00){\circle*{2.00}}
\put(24.33,26.00){\circle*{2.00}}
\put(24.33,36.00){\circle*{2.00}}
\put(24.33,1.00){\circle*{2.00}}
\put(39.33,6.00){\circle*{2.00}}
\put(9.33,6.00){\circle*{2.00}}
\emline{4.33}{21.00}{1}{24.33}{36.00}{2}
\emline{24.33}{36.00}{3}{44.33}{21.00}{4}
\emline{44.33}{21.00}{5}{24.33}{26.00}{6}
\emline{24.33}{26.00}{7}{4.33}{21.00}{8}
\emline{4.33}{21.00}{9}{9.33}{6.00}{10}
\emline{9.33}{6.00}{11}{24.33}{1.00}{12}
\emline{24.33}{1.00}{13}{39.33}{6.00}{14}
\emline{39.33}{6.00}{15}{44.33}{21.00}{16}
\put(1.00,21.00){\makebox(0,0)[cc]{\underline{\bf 1}2}}
\put(24.66,23.33){\makebox(0,0)[cc]{\underline{\bf 2}3}}
\put(28.33,35.67){\makebox(0,0)[cc]{1\underline{\bf 2}}}
\put(47.66,21.00){\makebox(0,0)[cc]{1\underline{\bf 3}}}
\put(36.00,7.33){\makebox(0,0)[cc]{\underline{\bf 1}3}}
\put(24.33,3.67){\makebox(0,0)[cc]{2\underline{\bf 3}}}
\put(12.33,7.33){\makebox(0,0)[cc]{1\underline{\bf 2}}}
\end{picture}
\caption{The graph $\theta_{2,2,4}$}
\label{theta224} %
\end{figure}
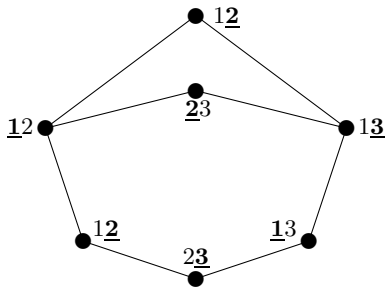

The foregoing two theorems are major steps in the proof of
Theorem~\ref{main}.
Before we proceed we must recall the definition of a $\theta$-graph.
If $p$, $q$, and $r$ are positive integers and at most one of them
equals~$1$, by $\theta_{p,q,r}$ we mean a graph which consists of
three internally disjoint paths of length $p$, $q$, and $r$ which
have the same endpoints. For example the graph $\theta_{2,2,4}$ is
shown in Figure~\ref{theta224}.

\begin{thm}
Suppose that $G$ is a $2$-connected graph, $t=\max\{3,\Chi(G)\}$,
and $G$ is not \uttlc. Then $G$ is either a complete or a
triangle-free graph.
\label{tfree}
\end{thm}

\begin{proof}{ %
Let $G$ be a graph which is not \uttlc\ for $t=\max\{3,\Chi(G)\}$,
and suppose that $G$ contains a triangle. For every pair of
independent vertices of $G$, say $u$ and $v$, which take on
different colors in each $t$-coloring of $G$, we add the edge $u
v$, to obtain a graph $G^*$. By Lemma~\ref{adduv}, $G^*$ is not a
\uttlc\ graph. If $G^*$ is not a complete graph, since it is
$2$-connected and contains a triangle, it must have an induced
$\theta_{1,2,r}$ subgraph, say $H$ (to see this, consider a
maximum clique in $G^*$ and a minimum path outside it which
joins
two vertices of this clique). Suppose that $x, y$, and $z$ are the
vertices of a triangle in $H$, and
$y=v_0,v_1,\ldots,v_{r-1},v_r=z$ is a path of length $r$ in $H$
not passing through $x$. Consider a $t$-coloring $c$ of $G^*$ in
which $x$ and $v_{r-1}$ take on the same color. We define a
$2$-list assignment $L$ to $H$ as follows.
$$L(x)=L(z)=\{c(x), c(z)\}, L(y)=\{c(x), c(y)\},$$
$$L(v_i)=\{c(v_i), c(v_{i-1})\} ;\ \ \forall\ 1\le i\le r-1.$$ In
each $L$-coloring of $H$ one of the vertices $x$ and $z$ must take
on the color $c(x)$ and the other takes on the color $c(z)$. So $y$
must take on the color $c(y)$ and one can see by induction that
each
$v_i$ must take on the color $c(v_i)$, and finally $x$ must take on
the color $c(x)$.
Now since $G^*$ is connected, as in the proof of
Lemma~\ref{one1}, one can extend $L$ to a $2$-list assignment to
$G^*$ such that $c$ is the only $L$-coloring of
$G^*$. This contradiction implies that $G^*$ is a complete graph,
and this means that $G$ has chromatic number $n(G)$, so $G$
must
be a complete graph.
}\end{proof}

\begin{thm}
\label{cyclechord}Let $G$ be a triangle-free $2$-connected graph
which contains a cycle with a chord and $t=\max\{3,\Chi(G)\}$. Then
$G$ is \uttlc\ if and only if it is not a complete bipartite graph.
\end{thm}

\begin{proof}{ %
By Theorem~\ref{mahu2lc}, a complete bipartite graph is not \utlc.
So if $G$ is \uttlc, it is not a complete bipartite graph. For the
converse, let $G$ be a graph which is not \uttlc\ where
$t=\max\{3, \Chi(G)\}$, and suppose that $G$ contains a cycle with
a chord. For every pair of independent vertices of $G$, say $u$
and $v$, which take on different colors in each $t$-coloring of
$G$, we add the edge $u v$, to obtain a graph $G^*$. By
Lemma~\ref{adduv}, $G^*$ is not a \uttlc\ graph. If $G^*$ contains
a triangle, By Theorem~\ref{tfree}, $G^*$ and so $G$ must be
complete graphs which contradicts the hypothesis. So suppose that
$G^*$ does not contain a triangle.

Consider a cycle $v_1 v_2\ldots v_p v_1$ with a chord $v_1v_\ell$,
and suppose $H$ to be the graph $G^*[v_1, v_2, \ldots, v_p]$. If
$v_pv_{\ell-1}\not\in E(H)$, there exists a $t$-coloring $c$ of
$G^*$, such that $c(v_p)=c(v_{\ell-1})$. Assign the list
$L(v_i)=\{c(v_i), c(v_{i-1})\}$ to each $v_i$, where $1\le i\le p$
and $v_0=v_p$. Consider an $L$-coloring $c'$ for $H$. Starting
from $v_1$ and considering each of two possible colors for it,
we conclude that $c'(v_\ell)=c(v_\ell)$. So
for each $1\le i\le p$ we have $c'(v_i)=c(v_i)$. This means that
$H$ is a \uttlc\ graph, and similar to the proof of Lemma~\ref{one1},
$G^*$ is a \uttlc\ graph, a contradiction. So $v_pv_{\ell-1}\in E(H)$
and similarly $v_2v_{\ell+1}\in E(H)$. Now consider the cycle
$v_1v_2v_{\ell+1}v_\ell v_{\ell-1}v_pv_1$ with chord $v_1v_\ell$. By
a similar
argument, $v_pv_{\ell+1}$ and $v_2v_{\ell-1}$ are in $E(H)$ and so
the graph $G^*[v_1v_2v_{\ell+1}v_lv_{\ell-1}v_p]$ is a $K_{3,3}$.

Suppose that $K$ is a maximal complete bipartite subgraph of $G^*$
containing the $K_{3,3}$ determined above.
Since $G$ is triangle-free, $K$ is an induced subgraph of $G$.
If $V(G)\setminus
V(K)\not=\emptyset$, consider a vertex $v\in V(G)\setminus V(K)$
which is adjacent to a vertex $w_1$ of $K$. By $2$-connectivity
of $G^*$, there exists a path $v u_1\ldots u_r w_2$ in which
$w_2\in V(K)$ and $u_i\not\in V(K)$ for each $0\le i\le r$. If
$w_1$ and $w_2$ are in the same part of $K$, since each part of
$K$ has at least $3$ vertices, there exists a vertex $w_3$ other
than $w_1$ and $w_2$ in the same part of $K$ as $w_1$ and $w_2$, and
vertices $w'_1$ and $w'_2$ in the other part of $K$. Considering
the cycle $v u_1\ldots u_rw_2w'_2w_3w'_1w_1v$ with chord $w_1w'_2$,
by a similar argument as in the previous paragraph,
it is implied that $v$ is adjacent to $w_3$. So $v$ is adjacent
to all the vertices of $K$ which are in the same part of $K$ as $w_1$,
except possibly to $w_2$, but in fact $v$ is adjacent to $w_2$,
since we can now consider $w_3$ in place of $w_2$ and do the
same as above. This contradicts the maximality of $K$.
On the other hand if $w_1$ and $w_2$ are in different parts of
$K$, a similar argument yields a contradiction.

We showed that $G^*=K$ and it is remained only to show that
$G=G^*$. If $x y$ is an edge in $G^*$ which is not present in $G$,
using the fact that $G$ is bipartite, one can easily obtain a
$t$-coloring $(t=3)$ of $G$ in which $x$ and $y$ take on the
same color, a contradiction.
}\end{proof}

At this point we will consider graphs that do not satisfy the
conditions of Theorem~\ref{cyclechord}, namely $2$-connected graphs
in which every cycle is induced.
The following lemma helps us to treat such graphs.
\begin{lem}
\label{deg2}A $2$-connected graph in which each cycle is
chordless, has at least a vertex of degree $2$.
\end{lem}

\begin{proof}{ %
It is a well-known theorem of H.~Whitney~\cite{whitney} that a
graph is $2$-connected, if and only if it admits an ear
decomposition (For a description of ear decomposition see
Theorem~4.2.7 in~\cite{west}). In the case of present lemma,
since the graph is chordless, each ear
is a path of length at least $2$, so the last ear contains a
vertex of degree $2$.
}\end{proof}

If $G$ is a graph and $v$ a vertex of $G$, we define $G_v$ to be a
graph obtained by identifying $v$ and all of its neighbors to a
single vertex $[v]$.
\begin{lem}
If $v$ is a vertex of degree $2$ in a graph $G$, and $G_v$ is
\uttlc\ for some $t$, then $G$ is also \uttlc.
\label{Gv}
\end{lem}

\begin{proof}{
Suppose that $v_1$ and $v_2$ are the neighbors of $v$ in
$G$. If $L$ is a $(2,t)$-list assignment to $G_v$
such that $G_v$ has a unique $L$-coloring, one can assign $L(w)$
to each vertex $w$ of the graph $G$ except $v$, $v_1$, and $v_2$, and
$L([v])$ to these three vertices, to obtain a $(2,t)$-list assignment
to $G$ from which $G$ has a unique list coloring.
}\end{proof}

The following lemma gives us a family of uniquely $(2,3)$-list colorable
graphs, which we will use in the proof of our main result.
\begin{lem}
Aside from $\theta_{2,2,2}=K_{2,3}$, each graph $\theta_{p,q,r}$ is
uniquely $(2,3)$-list colorable.
\label{theta}
\end{lem}

\begin{proof}{
Suppose that $G=\theta_{p,q,r}$ is a counterexample with minimum
number of vertices, and $u$ and $v$ are the two vertices of $G$
with degree $3$. If one of $p$, $q$, and $r$ is $1$, then $G$ is
a cycle with a chord and we have nothing to prove. Otherwise
suppose that one of the numbers $p$, $q$, and $r$, say $p$
is odd, and there exists a vertex $w$ on a path with length $p$
between $u$ and $v$. Then by Lemma~\ref{Gv} the graph $G_w$ is not
a uniquely $(2,3)$-list colorable graph, a contradiction. Hence
$p=1$ and we yield to the previous case.

So assume that $p, q$, and $r$ are all even numbers. By the
hypothesis at least one of $p$, $q$, and $r$, say $r$, is greater
than $2$. If either $p>2$, $q>2$, or $r>4$, by use of
Lemma~\ref{Gv} we obtain a smaller counterexample to the
statement, which is impossible by minimality of $G$, so
$G=\theta_{2,2,4}$. In Figure~\ref{theta224} there is given a
$(2,3)$-list assignment to $\theta_{2,2,4}$ which
induces a unique list coloring. This shows that $G$ is a
uniquely $(2,3)$-list colorable graph, which contradicts the
fact that $G$ is a counterexample to the statement.
}\end{proof}

Now we can prove the main result.
\begin{thm}\label{main}{\rm (MAIN)}
A graph $G$ is \utlc\ if and only if it is \uttlc, where
$t=\max\{3,\Chi(G)\}$.
\end{thm}

\begin{prooff}{ %
By definition, if $G$ is \uttlc\ for some $t$, it is \utlc. So we
must only prove that every \utlc\ graph $G$ is \uttlc\ for
$t=\max\{3,\Chi(G)\}$. Suppose that $G$ is a counterexample to the
statement with minimum number of vertices. By Theorem~\ref{block},
$G$ is $2$-connected, by Theorem~\ref{tfree}, it is
triangle-free (by Theorem~\ref{mahu2lc} it can not be a complete
graph), and by Theorem~\ref{cyclechord}, it does not have
a cycle with a chord, so Lemma~\ref{deg2} implies that $G$ has a
vertex~$v$ with exactly two neighbors $v_1$ and $v_2$.

Consider the graph $H=G\setminus v$ and note that since $\deg
v=2$, we have $\max\{3,\Chi(H)\}=\max\{3,\Chi(G)\}$. So if $H$ is
\utlc, by minimality of $G$, the graph $H$ must be \uttlc, and
since $t\ge 3$ and $\deg v=2$, we conclude that $G$ is \uttlc, a
contradiction. Therefore $H$ is not a \utlc\ graph and because it
is a triangle-free graph, by Theorem~\ref{mahu2lc} every block of
$H$ is either a cycle of length at least four or a complete bipartite
graph. This shows that $t=3$.

We will show by case analysis that $G$ has an induced subgraph
$G'$ which is isomorphic to some $\theta_{p,q,r}\not=\theta_{2,2,2}$
(except in the case (i.2)).
The graph $G'$ is \uttlc\ by Lemma~\ref{theta}. Now a $(2,3)$-list
assignment to $G'$ with a unique list coloring can simply be
extended to the whole of $G$. This completes the proof.
To show the existence of $G'$ we consider two cases.
\begin{itemize}
\item[(i)] The graph $H$ is $2$-connected.
So $H$ is either a $K_2$, a cycle, or a complete bipartite graph
with at least two vertices in each part. If $H=K_2$ then $G=K_3$,
a contradiction.
\begin{itemize}
\item[(i.1)] If $H$ is a cycle, $G$ is a $\theta$-graph and $G'=G$.
Note that since $G$ is not uniquely $2$-list colorable, $G'=G$ is
not isomorphic to $\theta_{2,2,2}$.
\item[(i.2)] If $H$ is a complete bipartite graph, since $G$ is
triangle-free, $v_1$ and $v_2$ are in the same part in $H$. Now
there must exist at least one other vertex $v_3$ in that part
--otherwise $G$ will be a complete bipartite graph. Suppose that
$u_1$ and $u_2$ are two vertices in the other part of $H$. The
graph $G'$ induced from $G$ on $\{v,v_1,v_2,v_3,u_1,u_2\}$ is a
uniquely $(2,3)$-list colorable with the list assignment $L$
as follows:
$L(v)=\{1,2\}$,
$L(v_1)=\{1,3\}$,
$L(v_2)=\{1,2\}$,
$L(v_3)=\{2,3\}$,
$L(u_1)=\{2,3\}$,
$L(u_2)=\{1,3\}$.
\end{itemize}
\item[(ii)] The graph $H$ is not $2$-connected. Since $G$ is
$2$-connected $H$ has exactly two end-blocks each of them contains
one of $v_1$ and $v_2$.

If all of the blocks of $H$ are isomorphic to $K_2$, then $G$ is a
cycle which is impossible. So $H$ has a block $B$ with at least
three vertices. Since $B$ is a cycle or a complete bipartite graph
with at least two vertices in each part, it has an induced cycle
$C$ which shares a vertex with at least two other blocks. Since
$G$ is $2$-connected these two vertices must be connected by a
path disjoint from $B$. Suppose that $P$ is such a path with
minimum length. The graph $G'=C\cup P$ is the required $\theta$-graph.
\hfill{$\blacksquare$}
\end{itemize}
}\end{prooff}

\section{Concluding remarks}
We begin with a definition which is a natural consequence of the
aforementioned results.
\begin{dfn}
For a graph $G$ and a positive integer $k$, we define $\Chi_u(G,k)$
to be the minimum number $t$, such that $G$ is a \uktlc\ graph,
and zero if $G$ is not a \uklc\ graph. The {\deffnt uniquely list
chromatic number} of a graph $G$, denoted by $\Chi_u(G)$, is
defined to be $\max_{k\ge 1}\Chi_u(G,k)$.
\end{dfn}

In fact Theorem~\ref{main} states that for every uniquely $2$-list
colorable graph~$G$,
$\Chi_u(G,2)=\max\{3, \Chi(G)\}$ and by Brooks' theorem and the
fact that for every \utlc\ graph $G$, $\Delta(G)\ge 3$, we have
shown that $\Chi_u(G,2)\le\Delta(G)+1$. This seems to remain true
if we substitute $2$ by any positive integer $k$.

\begin{conj}
For every graph $G$ we have $\Chi_u(G)\le\Delta(G)+1$, and equality
holds if and only if $G$ is either a complete graph or an odd cycle.
\end{conj}

The above conjecture implies the well-known Brooks' theorem, since
for every graph $G$ we have $\Chi_u(G,1)=\Chi(G)$, and so
$\Chi(G)\le\Chi_u(G)$. Hence the above conjecture implies that
$\Chi(G)\le\Delta(G)+1$.
On the other hand if $\Chi(G)=\Delta(G)+1$, we will have
$\Chi_u(G)=\Delta(G)+1$ and the conjecture above implies that $G$
is either a complete graph or an odd cycle.
\section*{Acknowledgements}
The authors are grateful to Professor E.\,S.~Mahmoodian for his
comments and support. We also thank the anonymous referees for their
inquiry and useful comments.

\addtolength{\baselineskip}{-1mm}
\newcommand{\noopsort}[1]{} \newcommand{\printfirst}[2]{#1}
  \newcommand{\singleletter}[1]{#1} \newcommand{\switchargs}[2]{#2#1}

\end{document}